\documentclass{amsart}
\usepackage{Everything}

\title{Universal countably chromatic graph}
\author{Siiri Kivimäki}
\thanks{This project has received funding from the European Research Council (ERC) under the
European Union’s Horizon 2020 research and innovation programme (grant agreement No
101020762) and Magnus Ehrnrooth Foundation.}

\date{\today}

\begin{document}

\maketitle

\begin{abstract}
    We show that the existence of a universal countably chromatic graph of size $\aleph_1$ together with the failure of continuum hypothesis is consistent. The proof is a forcing iteration of strongly proper ccc posets. The construction works for any uncountable successor cardinal $\kappa^+$, where $\kappa$ is regular.
\end{abstract}

\tableofcontents

\section{Introduction}

\noindent A graph is \textbf{countably chromatic} if there exists a coloring function on the set of vertices into $\N$ that assigns different colors to adjacent vertices. In this paper we consider the existence of a universal object among countably chromatic graphs. An object in a class equipped with embeddings is called \textbf{universal} if any other object in the class embeds into it. Here the objects are countably chromatic graphs of size $\aleph_1$ with a specified coloring; they are triples $G=(G,E,c)$, where $(G,E)$ is an undirected and loopless graph and $c:G\to\omega$ is a function such that for all vertices $a,b\in G$,
\[
c(a)=c(b)\implies\neg E(a,b).
\]
We refer to the function $c$ by \textbf{coloring} of $G$. An \textbf{embedding} between colored graphs $(G,E_G,c_G)$ and $(H,E_H,c_H)$ is an injective function $f:G\to H$ preserving the edge relation and the coloring:
\begin{enumerate}
    \item $E_G(a,b)$ if and only if $E_H(f(a),f(b))$,
    \item $c_G(a)=c_H(f(a))$.
\end{enumerate}
We refer to this class of objects and embeddings by \textbf{countably chromatic graphs of size $\bs{\aleph_1}$}. The main theorem of this paper is the following.

\begin{manualtheorem}{\ref{thm}}
    It is consistent that the class of countaby chromatic graphs of size $\aleph_1$ has a universal object and $2^{\omega}=\aleph_2$.
\end{manualtheorem}

The history of the universal countably chromatic graph of size $\aleph_1$ begins with the Rado graph, or random graph, which is a countable graph that is universal for the class of countable graphs. The first natural question was whether there exists an analogue of such graph in size $\aleph_1$. It turned out that the existence of a graph of size $\aleph_1$ that is both universal and \textit{strongly homogeneous}\footnote{every small partial elementary map lifts to an automorphism} is in fact equivalent to the continuum hypothesis (See Theorem 4.7 on page 476 in \cite{shelah1990classification}). Shelah \cite{shelah1990universal} then showed that the existence of a graph of size $\aleph_1$ that is universal but not necessarily strongly homogeneous is consistent with the failure of the continuum hypothesis. Also the non-existence is: after adding $\aleph_2$ many Cohen reals, there is no universal graph of size $\aleph_1$ and the continuum hypothesis fails. For a proof, see appendix of \cite{kojman1992nonexistence}. 

The purpose of this paper is to show that this picture generalizes to countably chromatic graphs. The continuum hypothesis implies the existence of a universal countably chromatic graph of size $\aleph_1$. The proof is similar as in the case of graphs, despite the fact the class of countably chromatic is not elementary (first-order definable). The trick is to move to a many-sorted signature with one sort for each color. By the continuum hypothesis, there is a saturated model of size $\aleph_1$ for this many-sorted theory, which is universal.

Likewise, it holds also that after adding $\aleph_2$ many Cohen reals there is no universal countably chromatic graph of size $\aleph_1$ and the proof is a straightforward modification of the proof in the appendix of \cite{kojman1992nonexistence}.

The main theorem (Theorem \ref{thm}) is proved by forcing. We start with a model of $\GCH$ and perform a finite support iteration of length $\omega_2$ of strongly proper ccc posets. The structure of the construction is inspired by \cite{ben2023aronszajn} and \cite{kivimaki2025universal} and the definition of the poset is inspired by \cite{shelah1990universal}. The same proof works for higher successor cardinals: for any cardinal $\kappa$, there is a poset that forces the existence of a graph of size $\kappa^+$ of chromatic number $\kappa$ that contains an isomorphic copy of each such graph, together with $2^\kappa>\kappa^+$.

\section{Strategy and preliminaries}

\noindent We assume $\GCH$ in the ground model.
The poset $\P$ will be a direct limit of a finite support iteration $(\P_\delta:\delta\leq\omega_2)$, where
\begin{enumerate}
    \item $\P_1$ builds a generic countably chromatic graph structure $\dot \GG$ on $\omega_1$,
    \item the tail of the iteration creates embeddings
    \[
    \dot f_\delta: \dot H_\delta\to \dot\GG,
    \]
    where $\dot H_\delta$ is a $\P_\delta$-name for a countably chromatic graph on $\omega_1$ chosen by a suitable bookkeeping function, for every $\delta<\omega_2$.
\end{enumerate}

\noindent The iteration will be defined in such a way that at each stage $\delta<\omega_2$, we will fix a sequence of countable models $(M^\delta_\alpha:\alpha\in\EE_\delta)$ indexed by a club $\EE_\delta\subseteq\omega_1$, and it will be made sure that the top condition $1_\P\in\P_\delta$ is strongly $(\P_\delta,M^\delta_\alpha)$-generic, for every $\alpha\in\EE_\delta$. The role of these models is similar to side conditions, with the exception that they do not explicitly appear in the poset. It will follow that the final poset $\P_{\omega_2}$  has ccc. 

\vv 

\noindent We recall the relevant definitions related to strong properness. If $\P$ is a poset and $p,q\in\P$ are conditions, we denote $p||q$ if $p$ and $q$ are compatible.

\begin{de}
    Let $\P$ be a poset and let $M$ be a set.
    \begin{enumerate}
        \item A condition $r\in\P\cap M$ is a \textbf{residue} of $p$ into $M$ if
        \[
        \forall w\in\P\cap M(w\leq r\to w\,||\,p).
        \]
        \item A condition $p\in\P$ is \textbf{strongly $\bs{(\P,M)}$-generic} if every $q\leq p$ has a residue into $M$.
        \item The poset $\P$ is \textbf{strongly proper with respect to $\bs{M}$} if for every $p\in\P\cap M$ there is $q\leq p$ that is strongly $(\P,M)$-generic.
    \end{enumerate}
    Furthermore, $\P$ is \textbf{strongly proper} if it is strongly proper with respect to club many $M\in[H(\theta)]^\omega$ for any large enough regular cardinal $\theta$.
\end{de}

We say that a model $M$ is \textbf{suitable} for a poset $\P$ if it is an elementary submodel of some $(H(\theta),\in)$, where $\theta$ is a large enough regular cardinal, and $\P\in M$.

\begin{lem}\label{lemstronglygencha} Let $\P$ be a poset and let $M$ be suitable for $\P$. The following are equivalent for a condition $p\in\P$:
    \begin{enumerate}
        \item $p$ is strongly $(\P,M)$-generic,
        \item the set of conditions that have a residue into $M$ is dense below $p$,
        \item $p\Vdash``\check G\cap M$ is $V$-generic on $\P\cap M$".
    \end{enumerate}
\end{lem}

\begin{lem}\label{lemccc} Let $\P$ be a poset. If there are club many countable $M\elem H(\theta)$ such that every condition has a residue into $M$, then $\P$ is strongly proper and has ccc.
\end{lem}
\begin{proof}
    If a condition $p$ is strongly $(\P,M)$-generic, then it is $(\P,M)$-generic, and being ccc is equivalent to the statement that every condition is $(\P,M)$-generic for club many countable $M\elem H(\theta)$.
\end{proof}

\begin{lem}
    Strongly proper posets preserve $\omega_1$ and add reals.
\end{lem}

\noindent More on strongly proper forcing can be found for example at \cite{gilton2017side}. Our notation is standard and follows \cite{jech2003set}.

\section{Poset}

\noindent The first step is to force a generic countably chromatic graph of size $\aleph_1$. This is done with finite conditions, with a poset equivalent to the poset $\Add(\omega,\omega_1)$ that adds $\omega_1$ may Cohen reals. We assume that $\GCH$ holds in the ground model.

\begin{de}[$\P_1$]\label{de:firstposet}
    The poset $\P_1$ consists of triples
    \[
    p=(u^p,E^p,c^p),
    \]
    where
    \begin{enumerate}
        \item $u^p\subseteq\omega_1$ is a finite set,
        \item $E^p:[u^p]^2\to 2$ is a partial function,
        \item $c^p:u^p\to\omega$ is a function such that if $c^p(s)=c^p(t)$ and $\{s,t\}\in\dom(E^p)$, then $E^p(\{s,t\})=0$.
    \end{enumerate}
    The ordering is by $q\leq p$ iff $u^q\supseteq u^p$, $E^q\supseteq E^p$ and $c^q\supseteq c^p$.
\end{de}

\noindent A generic filter $G\subseteq\P_1$ gives a graph on $\omega_1$ by letting
\[
E^G(s,t):\iff\exists p\in G\spa E^p(\{s,t\})=1,
\]
and the countable chromaticity of this graph is witnessed by the coloring
\[
c^G:=\bigcup_{p\in G}c^p.
\]

\begin{notation} Let $\dot \GG$ be a $\P_1$-name for the graph with a coloring $(\omega_1,E^{\check G},c^{\check G})$.
\end{notation}

The poset $\P_1$ is equivalent to the poset that adds a subset of $\omega_1$ with finite conditions, or equivalently to the poset $\Add(\omega,\omega_1)$ that adds $\omega_1$ many Cohen reals.
The goal is to define a poset $\P_{\omega_2}$ which creates an embedding 
\[
\dot f_\delta:\dot H_\delta\to \dot \GG
\]
for each $\delta<\omega_2$, where $\dot H_\delta$ is a countably chromatic graph on $\omega_1$ chosen by a suitable bookkeeping function. The first poset will be $\P_1$. Being countably chromatic is preserved by any forcing extension, so when forcing the embeddings, we only need to worry about the preservation of $\aleph_1$.

For each $\delta<\omega_2$, we define models $(M^\delta_\alpha:\alpha<\omega_1)$ which will implicitly play the role of side-conditions, even if they will not explicitly appear as part of the conditions in $\P_\delta$. For each poset $\P_\delta$, there will be a club $\EE_\delta\subseteq\omega_1$ such that $\P_\delta$ is strongly proper with respect to each model $M^\delta_\alpha$ for $\alpha\in\EE_\delta$. These clubs $\EE_\delta$ will get ``faster" in the sense that whenever $\alpha\in\EE_\delta$ and $\gamma\in\delta\cap M^\delta_\alpha$, then $\alpha$ is a limit point of $\EE_\gamma$.

\begin{notation} We fix a regular cardinal $\theta>\aleph_2$ and a well-ordering $<_\theta$ of $H(\theta)$.
    We also fix a bookkeeping function \[
    \omega_2\to H(\omega_2),\quad \gamma\mapsto \dot H_\gamma
    \] 
    such that whenever the poset $\P_\gamma$ is defined, then $\dot H_\gamma$ is a $\P_\gamma$-name for a countably chromatic graph on $\omega_1$. We assume that for every $\gamma<\omega_2$ and for every $\P_\gamma$-name for a countably chromatic graph on $\omega_1$ there is a name for an isomorphic copy of it in the range of the bookkeeping function. We assume that the bookkeeping function is a minimal such function with respect to $<_\theta$.
\end{notation}

The posets $\P_\delta$ and clubs $\EE_\delta\subseteq\omega_1$ are defined by recursion on $\delta<\omega_2$. The first poset $\P_0$ will be the trivial poset $\{\emptyset\}$. The poset $\P_1$ was already defined in \ref{de:firstposet}.

\vv

\noindent For a set $X\subseteq\omega_1$, the set $\lim X$ is the set of limit points of $X$.

\begin{de}
    Let $\delta<\kappa^+$. Whenever the poset $\P_\delta$ and the sets $\bar{\EE}_\delta:=(\EE_\gamma)_{\gamma<\delta}$ are defined, we let:
    \begin{enumerate}
        \item For every $\alpha\leq\kappa$, let \[M^\delta_\alpha:=\text{Skolem hull of }\alpha\cup\{\delta,\P_\delta,\bar{\EE}_\delta\}\text{ in }(H(\theta),\in,<_\theta),\]
    \item Let
        \begin{align*}
            \EE_\delta:=\{\alpha<\kappa:\kappa\cap M^\delta_\alpha=\alpha\text{ and }\alpha\in\bigcap_{\gamma\in\delta\cap M^\delta_\alpha}\lim\EE_\gamma\}.
        \end{align*}
        Moreover, for $\alpha<\omega_1$, let
    \begin{align*}
        &\rho_\delta(\alpha):=\text{ the least element in }\EE_\delta\text{ strictly above }\alpha,\\
        &\lambda_\delta(\alpha):=\text{ the least limit point of }\EE_\delta\text{ strictly above }\rho_\delta(\alpha).
    \end{align*}
    \end{enumerate}
\end{de}

\noindent Then each $\EE_\delta$ is a club, by normality of the club filter. For every $\alpha<\kappa$, we have
\[
\alpha<\rho_\delta(\alpha)<\lambda_\delta(\alpha).
\]
When mapping the $\delta$-th graph $\dot H_\delta$ on $\omega_1$ into $\dot\GG$, we will use finite partial functions $f:\omega_1\to\omega$ such that each $\alpha\in\dom(f)$ is mapped into the interval $f(\alpha)\in[\rho_\delta(\alpha),\lambda_\delta(\alpha))$. The ordinal $\rho_\delta(\alpha)$ is used as a ``rank" that separates each vertex $\alpha$ in $\dot H_\delta$ from its image $f(\alpha)$.

For a technical reason that appears in the proof of strong properness, we need to label each node in $\dot\GG$ with a countable ordinal. This label will be used to read off the $\rho_\delta(\alpha)$ from the preimage $\alpha$ of a node $\beta\in\dot\GG$.

\begin{notation}
    We fix a function $\lab:\omega_1\to\omega_1$ such that
    \begin{enumerate}
        \item $\lab(\alpha)<\alpha$ for every non-zero $\alpha<\omega_1$,
        \item for every $\alpha<\omega_1$ there are unboundedly many $\beta>\alpha$ with $\lab(\beta)=\alpha$.
    \end{enumerate}
    We call $l(\alpha)$ the \textbf{label} of $\alpha$.
\end{notation}

We are ready to define the embedding posets. The definition is recursive. In addition to the separation of each vertex $\alpha$ in the $\delta$-th graph $\dot H_\delta$ from its image in $\dot\GG$ using the ordinal $\rho_\delta(\alpha)$, we need to be able to read off $\rho_\delta(\alpha)$ from the image of $\alpha$. We will make sure to embed $\alpha$ to some $\beta$ whose label $\lab(\beta)$ is the ordinal $\rho_\delta(\alpha)$.

The clubs $\EE_\gamma$, $\gamma<\delta$ will be used in the definition of the poset $\P_\delta$, and the club $\EE_\delta$ is defined once the poset $\P_\delta$ is defined.

\begin{de}[$\P_\delta$]
    Let $\delta\leq\omega_2$. Assume that $\P_\gamma$ and $\EE_\gamma$ are defined for $\gamma<\delta$. The conditions in $\P_\delta$ are functions $p:\delta\to H(\omega_2)$ satisfying the following:
    \begin{enumerate}
        \item $p(0)=(u^p,E^p,c^p)\in\P_1$,
        \item For non-zero $\gamma<\delta$, $p(\gamma)$ is a pair 
        \[
        p(\gamma)=(f^p_\gamma,S^p_\gamma),
        \]
        where
        \begin{enumerate}
            \item $f^p_\gamma:\omega_1\to u^p$ is a finite partial injective function,
            \item $f^p_\gamma(\alpha)\in[\rho_\gamma(\alpha),\lambda_\gamma(\alpha))$ for all $\alpha\in\dom(f^p_\gamma)$,
            \item the label of $f^p_\gamma(\alpha)$ satisfies $\lab(f^p_\gamma(\alpha))=\rho_\gamma(\alpha)$,
            \item $p\rest\gamma$ decides the $\dot H_\gamma$-edges and colors in $\dom(f^p_\gamma)$,
            \item $p\rest\gamma\Vdash ``f^p_\gamma:\dot H_\gamma\to\dot\GG$ is edge- and color-preserving",
            \item $S^p_\gamma\subseteq\omega_1$ is a finite set and $\ran(f^p_\gamma)\cap S^p_\gamma=\emptyset$,
        \end{enumerate}
        \item the \textbf{support} $\sp(p):=\{\gamma<\delta:p(\gamma)\neq(\emptyset,\emptyset)\}$ is finite.
    \end{enumerate}
    The ordering is defined by $q\leq p$ if:
    \begin{enumerate}
        \item $q(0)\leq p(0)$,
        \item $f^q_\gamma\supseteq f^p_\gamma$ and $S^q_\gamma\supseteq S^p_\gamma$ for every non-zero $\gamma<\delta$.
    \end{enumerate}
\end{de}

The posets $\P_\delta$ as well as the auxiliary sets $\EE_\delta$, $\delta\leq\omega_2$, have now been defined. We will show first a density claim about the embeddings, and then that each poset $\P_\delta$ has ccc. This will be done by showing that every condition $p\in\P_\delta$ is strongly $(\P_\delta, M^\delta_\alpha)$-generic, for $\delta<\omega_2$ and $\alpha\in\EE_\delta$. It will follow that each $\P_\delta$ is strongly proper and has ccc, and thus $\P_{\omega_2}$ will have ccc, being the direct limit of ccc posets.

The following is clear from definitions:

\begin{lem}
    If $\gamma\leq\delta\leq\omega_2$, $p\in\P_\delta$ and $q\in\P_\gamma$ extends $p\rest\delta$, then the concatenation
    \[
    q^\smallfrown p\rest[\gamma,\delta)
    \]
    is a condition in $\P_\delta$ that extends $p$.
\end{lem}

\noindent In particular, each poset $\P_\gamma$ embeds completely into $\P_\delta$, for $\gamma\leq\delta\leq\omega_2$.

\section{Density}

\noindent In this section we show that in $V^{\P_{\omega_2}}$, the graph $\dot\GG$ contains an isomorphic copy of each graph $\dot H_\delta$. In other words, we show that if $G\subseteq\P_{\omega_2}$ is a generic, then for each $\delta<\omega_2$, the map $f^G_\delta:=\bigcup_{p\in G}f^p_\delta$ satisfies $\dom(f^G_\delta)= H_\delta^G$. It suffices to show that for every $\delta<\omega_2$ and $\alpha<\omega_1$ the set $\{p\in\P_{\omega_2}:\alpha\in\dom(f^p_\delta)\}$ is dense.

\vv

The following is immediate:

\begin{lem} Let $G\subseteq\P_{\omega_2}$ be a generic filter.
\begin{enumerate}
    \item For every $\delta<\omega_2$, the function
    \[
    f^G_\delta:=\bigcup_{p\in G}f^p_\delta:\dot H_\delta^G\to \dot\GG^G
    \]
    is edge- and color-preserving.
    \item The function 
    \[
    c^G:=\bigcup_{p\in G}c^p:\dot\GG^G\to\omega
    \]
    is a coloring of $\dot\GG^G$, witnessing that $\dot\GG^G$ is countably chromatic.
\end{enumerate}
\end{lem}
\noindent We show that the generic embeddings are defined on the whole domain of $\dot H_\delta$ for each $\delta<\omega_2$.

\begin{lem}
    For every $p\in\P_{\omega_2}$, $\delta<\omega_2$ and every $\alpha<\omega_1$ there is $q\leq p$ such that $\alpha\in\dom(f^q_\gamma)$.
\end{lem}
\begin{proof}
    Up to extending $p\rest\delta$, we may assume that it decides the $\dot H_\delta$-edges and colors in the set $\dom(f^p_\delta)\cup\{\alpha\}$. The interval $[\rho_\gamma(\alpha),\lambda_\gamma(\alpha))$ is an infinite subset of $\GG$ and the set $u^p\cup S^p_\delta$ is finite, so we may find a node $\beta\in[\rho_\gamma(\alpha),\lambda_\gamma(\alpha))-(u^p\cup S^p_\delta)$ and extend $p(0)$ by assigning the color of $\alpha$ to  $\beta$ and connecting it to those nodes $f^p_\delta(\alpha')\in u^p$ such that $p\rest\delta\Vdash``E_{\dot H_\delta}(\alpha,\alpha')"$ and disconnect it from all other nodes in $u^p$. Define $q$ such that $q(0)$ is the resulting extension of $p(0)$ and
    \begin{align*}
        & q(\gamma)=\begin{cases}
            (f^p_\gamma\cup\{(\alpha,\beta)\},S^p_\gamma),\quad &\text{if }\gamma=\delta,\\
            p(\gamma), &\text{otherwise.}
        \end{cases}
    \end{align*}
    This $q$ is a condition, extends $p$ and has $\alpha$ in the domain of $f^q_\delta$.
\end{proof}

\section{Strong properness and ccc}

It remains to show that $\omega_1$ is preserved and that $2^\omega=\omega_2$ in $V^{\P_{\omega_2}}$. This is done by showing that each $\P_\delta$ is strongly proper with respect to the models 
\[
\{M^\delta_\alpha:\alpha\in\EE_\delta\}.
\]
In fact, we will show that the top condition $1_{\P_\delta}$ is strongly $(\P_\delta, M^\delta_\alpha)$-generic, for each $\alpha\in\EE_\delta$. This implies that each $\P_\delta$ has ccc (see \cite{shelah2017proper}). In order to show that $1_{\P_\delta}$ is strongly $(\P_\delta,M^\delta_\alpha)$-generic, it suffices to show that for densely many conditions $p$, the pointwise intersection of $p$ with $M^\delta_\alpha$, its \textit{trace}, is a residue of $p$ into $M^\delta_\alpha$.

\begin{de}
    Let $\delta<\omega_2$ and $\alpha\in\EE_\delta$. The \textbf{trace} of a condition $p\in\P_\delta$ into $M^\delta_\alpha$ is defined to be
    \begin{align*}
        [p]^\delta_\alpha:\delta\to M^\delta_\alpha,\quad p(\gamma):=\begin{cases}
            (u^p\cap M^\delta_\alpha,E^p\cap M^\delta_\alpha,c^p\cap M^\delta_\alpha)\quad &\text{if }\gamma=0,\\
            (f^p_\gamma\cap M^\delta_\alpha,S^p_\gamma\cap M^\delta_\alpha), &\text{if }\gamma\in\delta\cap M^\delta_\alpha-\{0\},\\
            \emptyset, &\text{if }\gamma\notin\delta\cap M^\delta_\alpha.
        \end{cases}
    \end{align*}
\end{de}

\noindent The trace might not even be a condition. It will be one whenever $p$ is \textit{super-nice}, in which case it will be a residue of $p$ as well.

\begin{rmk}\label{rmk:closed}
    It follows from the definition of the poset that if $\delta<\omega_2$, $\alpha\in\EE_\delta$, $p\in\P_\delta$ and $\gamma\in\delta\cap M^\delta_\alpha$, then the model $M^\delta_\alpha$ is closed for the function $f^p_\gamma$: if $\alpha'\in\dom(f^p_\gamma)\cap M^\delta_\alpha$, then $\alpha'<\alpha$. Since $\gamma\in M^\delta_\alpha$, the ordinal  $\alpha$ must be a limit point of $\EE_\gamma$. Thus $\lambda_\gamma(\alpha')\leq\alpha$. Now 
    \[
    f^p_\gamma(\alpha')\in[\rho_\gamma(\alpha'),\lambda_\gamma(\alpha'))\subseteq\alpha\subseteq M^\delta_\alpha.
    \]
So $f^p_\gamma(\alpha')\in M^\delta_\alpha$.
\end{rmk}

\begin{rmk}
     The sets $S^p_\gamma$ are used as follows. Suppose that $p\in\P_\delta$ is a condition, $\alpha\in\EE_\delta$, and we are showing that the trace $[p]^\delta_\alpha$ is a residue of $p$ into $M$. Let $\gamma<\delta$, $\gamma\in M^\delta_\alpha$. Suppose for simplicity that the graph $\dot H_\gamma$ is in $V$. By Remark \ref{rmk:closed} above the model $M^\delta_\alpha$ is closed under the function $f^p_\gamma$.
    \begin{itemize}
        \item Suppose that $\beta_1\in \omega_1\cap M^\delta_\alpha$, $\beta_2\in\omega_1-M^\delta_\alpha$ and $f^p_\gamma$ maps some node $\alpha_2\in M$ to $\beta_2$, but no node to $\beta_1$.
        \item Suppose that $w\in\P_\delta\cap M^\delta_\alpha$ extends $[p]^\delta_\alpha$ by sending some $\alpha_2$ to $\beta_2$. 
        \item If there is an edge $E_{H_\gamma}(\alpha_1,\alpha_2)$ but $p\Vdash \neg E^{\check G}(\beta_1,\beta_2)$, or vice versa, then $w$ cannot be compatible with $p$.
        \item We put into $S^p_\gamma$ every node in $M^\delta_\alpha-\ran(f^p_\gamma)$ of which $p$ decides that it is connected or disconnected from a vertex outside of $M$. There are finitely many of these and they all are in $u^p$. 
        \item In particular, the vertex $\beta_2$ is not in $M^\delta_\alpha$. If $p$ decides the value of $E^{\check G}(\beta_1,\beta_2)$, then $\beta_1$ is either already in the range of $f^p_\gamma$, or in the set of forbidden nodes $S^p_\gamma$. 
        \item So if $w\leq [p]^\delta_\alpha$, then $S^w_\gamma\supseteq S^p_\gamma\cap M$, so either $w$ agrees with $p$ on the node mapped to $\beta_1$, or else it cannot map anything on $\beta_1$.
    \end{itemize}
\end{rmk}

\vv 

We will show that every condition in $\P_\delta$ can be extended into a condition which is \textit{super-nice} with respect to $M^\delta_\alpha$, and that whenever $p$ is super-nice with respect to $M^\delta_\alpha$, then its trace is a residue of $p$.

\begin{de} The definition is by recursion on $\delta<\omega_2$. Let $\alpha\in\EE_\delta$. A condition $p\in\P_\delta$ is \textbf{super-nice with respect to $\bs{M^\delta_\alpha}$} if
\begin{enumerate}
    \item For every $\gamma\in\sp(p)\cap M^\delta_\alpha$ and $\beta\in\ran(f^p_\gamma)-M^\delta_\alpha$: if $\{\beta,\beta'\}\in\dom(E^p)$ for some $\beta'\in M^\delta_\alpha$, then $\beta'\in \ran(f^p_\gamma)\cup S^p_\gamma$,
    \item if $\gamma\in\delta\cap M^\delta_\alpha$ and $\alpha'\in\dom(f^p_\gamma)-M^\delta_\alpha$, then $p\rest\gamma$ is super-nice with respect to $M^\gamma_{\rho_\gamma(\alpha')}$.
\end{enumerate}
\end{de}

\begin{lem}
    Let $\delta<\omega_2$ and $\alpha\in\EE_\delta$. If $p\in\P_\delta$ is super-nice with respect to $M^\delta_\alpha$ and $\gamma\in\delta\cap M^\delta_\alpha$, then $p\rest\gamma$ is super-nice with respect to $M^\gamma_\alpha$.
\end{lem}

We give a technical but practical characteriztion of super-niceness in terms of finite ``paths". A \textbf{path in $\bs{p}$ from $\bs{M^\delta_\alpha}$ to $\bs{M^\gamma_\rho}$} is a finite sequence of pairs $(\gamma_0,\rho_0),\dots,(\gamma_n,\rho_n)$ satisfying:
\begin{enumerate}
    \item $\gamma_0=\delta$ and $\rho_0=\alpha$, and $\gamma_n=\gamma$ and $\rho_n=\rho$,
    \item $\gamma_{k+1}\in\sp(p\rest\gamma_k)\cap M^{\gamma_k}_{\rho_k}$,
    \item there is $\alpha'\in\dom(f^p_{\gamma_{k+1}})-\rho_k$ such that $\rho_{k+1}=\rho_{\gamma_{k+1}}(\alpha')$. 
\end{enumerate}

\begin{lem}
    Let $\delta<\omega_2$, $p\in\P_\delta$ and $\alpha\in\EE_\delta$. The following are equivalent:
    \begin{enumerate}
        \item $p$ is super-nice with respect to $M^\delta_\alpha$,
        \item For every path $(\gamma_0,\rho_0),\dots,(\gamma_n,\rho_n)$ in $p$ from $M^\delta_\alpha$ to $M^\gamma_\rho$,
        it holds that for all $\xi\in\sp(p\rest\gamma_n)\cap M^{\gamma_n}_{\beta_n}$: if $\beta\in \ran(f^p_\xi)-{\rho_n}$ and  $\alpha<\rho_n$ are such that $\{\alpha,\beta\}\in\dom(E^p)$, then $\alpha\in\ran(f^p_\xi)\cup S^p_\xi$.
    \end{enumerate}
\end{lem}

\begin{lem}
    Let $\delta<\omega_2$ and $\alpha\in\EE_\delta$. For every $p\in\P_\delta$ there is $q\leq p$ which is super-nice with respect to $M^\delta_\alpha$.
\end{lem}
\begin{proof}
    Let $p\in\P_\delta$. We extend $p$ by adding all nodes in $u^p-\ran(f^p_\gamma)$ into $S^p_\gamma$ for every $\gamma<\delta$. This is an overkill but causes no harm. Indeed, define $q\in\P_\delta$ by letting $q(0):=p(0)$ and for non-zero $\gamma<\delta$, let
    \[
    q(\gamma):=(f^p_\gamma,S^p_\gamma\cup(u^p-\ran(f^p_\gamma)).
    \]
    Then $q$ is super-nice with respect to $M^\delta_\alpha$. In fact, then $q$ is super-nice with respect to every structure.
\end{proof}

We are ready to show that if $p$ is super-nice with respect to $M^\delta_\alpha$, then its trace $[p]^\delta_\alpha$ is its residue into $M^\delta_\alpha$.

\begin{prop}\label{prop:sp}
    Let $\delta<\omega_2$ and $\alpha\in\EE_\delta$.
    Assume that $p\in\P_\delta$ is super-nice with respect to $M^\delta_\alpha$. Then
    \begin{enumerate}
        \item $[p]^\delta_\alpha\in\P_\delta\cap M^\delta_\alpha$,
        \item $\forall w\in\P_\delta\cap M^\delta_\alpha(w\leq [p]^\delta_\alpha\to w\,||\, p)$.
    \end{enumerate}
\end{prop}
\begin{proof} The proof is by induction on $\delta$.

\vv 

\noindent\textbf{Case $\bs{\delta=1}$:} 

\vv

\noindent It is clear that the trace $[p]^1_\alpha$ is a condition in $\P_1\cap M^1_\alpha$. If $w\in\P_1\cap M^1_\alpha$ extends the trace $[p]^\delta_\alpha$, then the pointwise union $w\cup p:=(u^w\cup u^p,E^w\cup E^p,c^w\cup c^p)$ is a common extension of $w$ and $p$, for all $p\in\P_1$ and $w\in\P_1\cap M^1_\alpha$, $\alpha<\omega_1$. Note that the edge function $[u^w\cup u^p]^2\to 2$ of a condition is only required to be a partial function. By taking a pointwise union, we do not introduce any new edges, so the function $c^w\cup c^p$ does not assign the same color to adjacent vertices.

\vv 

\noindent\textbf{Case $\bs{\delta}$ limit:} 

\vv 

\noindent First we show a claim about the structure of the models.

\begin{claim}\label{claim:limitcase}
    If $A$ is a finite subset of $\delta$ in $M^\delta_\alpha$, then there is $\gamma\in\delta\cap M^\delta_\alpha$ such that $A\subseteq\gamma\cap M^\gamma_\alpha$.
\end{claim}
\begin{proof}
    We show that if $\xi<\gamma<\delta$ are such that $\xi,\gamma\in M^\delta_\alpha$, then $\xi\in M^\gamma_\alpha$. To this end, fix a bijection $\psi:\omega_1\to\gamma$. (If $\gamma<\kappa$, then $\gamma<\alpha$, and there is nothing to prove.) We may assume that this bijection is $<_\theta$-minimal, so that it belongs to every elementary submodel of $(H_\theta,\in,<_\theta)$ containing $\gamma$. In particular $\psi\in M^\gamma_\alpha$. Since $\xi\in \gamma\cap M^\delta_\alpha$, there is $\beta<\alpha$ such that $\psi(\beta)=\xi$. Since $\beta,\psi\in M^\gamma_\alpha$, also $\xi=\psi(\beta)\in M^\gamma_\alpha$.
    
    Now the claim follows by letting $\gamma:=\textsf{max}(A)+1$.
\end{proof}

\noindent We finish the limit case.

\begin{itemize}
	\item To see that $[p]^\delta_\alpha$ is a condition: Using Claim \ref{claim:limitcase}, find $\gamma\in\delta\cap M^\delta_\alpha$ such that $\sp([p]^\delta_\alpha)\subseteq\gamma\cap M^\gamma_\alpha$. Then also $[p\rest\gamma]^\gamma_\alpha=[p]^\delta_\alpha\rest\gamma$. It follows that \[
[p]^\delta_\alpha={[p\rest\gamma]^\gamma_\alpha}^\smallfrown((\emptyset,\emptyset))_{\xi\in[\gamma,\delta)}.
\]
Since $[p\rest\gamma]^\gamma_\alpha$ is a condition by induction hypothesis, so must be $[p]^\delta_\alpha$.
	\item To see that $[p]^\delta_\alpha$ is a residue of $p$: Let $w\in\P_\delta\cap M^\delta_\alpha$ extend $[p]^\delta_\alpha$. Let $\gamma\in\delta\cap M^\delta_\alpha$ be such that $\sp(w)\subseteq\gamma\cap M^\delta_\alpha$. By Claim \ref{claim:limitcase} we have $\sp(w)\subseteq M^\gamma_\alpha$. By induction hypothesis there is $q_0\in\P_\gamma$ extending $w\rest\gamma$ and $p\rest\gamma$. Then $q_0$ can be extended into a condition $q\in \P_\delta$ extending $w$ and $p$ by letting
\[
q(\xi):=\begin{cases}
    q_0(\xi)\quad&\text{if }\xi<\gamma,\\
    p(\xi) &\text{if }\xi\in[\gamma,\delta).
\end{cases}
\]
\end{itemize}

\vv 

\noindent\textbf{Case $\bs{\delta+1}$:}

\vv

\noindent We argue first that $[p]^{\delta+1}_\alpha$ is a condition in $\P_{\delta+1}\cap M^{\delta+1}_\alpha$. Since $[p]^{\delta+1}_\alpha$ is a finite function $\delta+1\to M^{\delta+1}_\alpha$, it must be element of $M^{\delta+1}_\alpha$. The restriction $[p\rest\delta]^\delta_\alpha$ is a condition in $\P_\delta\cap M^\delta_\alpha$ and residue of $p\rest\delta$ into $M^\delta_\alpha$, by induction hypothesis and by the fact that $M^{\delta+1}_\alpha=M^\delta_\alpha$. Since
\[
[p]^{\delta+1}_\alpha={[p\rest\delta]^\delta_\alpha}^\smallfrown(f^{[p]^{\delta+1}_\alpha}_\delta,S^{[p]^{\delta+1}_\alpha}_\delta),
\]
it suffices to show that
\begin{enumerate}
    \item $f^p_\delta\rest\alpha$ is a function in $M^{\delta}_\alpha$,
    \item the condition $[p\rest\delta]^\delta_\alpha$ decides edges and colors in the domain of $f^{[p]^{\delta+1}_\alpha}_\delta$ and forces that it is edge- and color-preserving.
\end{enumerate}
The first item follows from Remark \ref{rmk:closed}. The second item follows using the induction hypothesis - since $[p\rest\delta]^\delta_\alpha$ is a residue of $p\rest\delta$ into $M^\delta_\alpha$, it forces the same facts about $f^{[p]^{\delta+1}_\alpha}_\delta$ as $p\rest\delta$ does. This is enough to see that $[p]^{\delta+1}_\alpha$ is a condition in $\P_{\delta+1}\cap M^{\delta+1}_\alpha$. There remains to show that it is a residue of $p$.

\vv

We then argue that $[p]^{\delta+1}_\alpha$ is a residue of $p$ into $M^{\delta+1}_\alpha$. Let $w\in\P_{\delta+1}\cap M^{\delta+1}_\alpha$ extend $[p]^{\delta+1}_\alpha$. We find a common extension $q\leq w,p$. 
We already have:
\begin{itemize}
    \item $w\rest\delta$ and $p\rest\delta$ are compatible, by induction hypothesis,
    \item $f^w_\delta\cup f^p_\delta$ is a finite injection,
    \item $(u^w\cup u^p,E^w\cup E^p,c^w\cup c^p)$ is a condition in $\P_1$,
    \item $(S^w_\delta\cup S^p_\delta)\cap\ran(f^w_\delta\cup f^p_\delta)=\emptyset$.
\end{itemize}
We need to find an extension $v\leq w\rest\delta,p\rest\delta$ which decides edges in the set $\dom(f^w_\delta)\cup\dom(f^p_\delta)$ such that $v(0)$ decides edges in $u^w\cup u^p$ \textit{correctly}: in such a manner that makes $f^w_\delta\cup f^p_\delta$ edge-preserving.
This condition $v$ is found in finitely many steps, climbing up the models $\{M^\delta_{\rho_\delta(\alpha')}:\alpha'\in\dom(f^p_\delta)-\alpha\}$. Each model $M^\delta_{\rho_\delta(\alpha')}$ separates the node $\alpha'$ fromits image  $f^p_\delta(\alpha')$ in the sense that \[
\alpha'\in M^\delta_{\rho_\delta(\alpha')}\quad\text{and}\quad f^p_\delta(\alpha')\notin M^\delta_{\rho_\delta(\alpha')}.
\] 
The point is to decide $\dot H_\delta$-edges of the set $\dom(f^w_\delta)\cup \{\alpha'\}$ inside the model $M^\delta_{\rho_\delta(\alpha')}$ and then extend $w(0)\cup p(0)$ accordingly outside of the model $M^\delta_{\rho_\delta(\alpha')}$, using the fact that $f^p_\delta(\alpha')\notin M^\delta_{\rho_\delta(\alpha')}$.

\vv 

\noindent Let first
\begin{align*}
    &f:=f^w_\delta\cup f^p_\delta.\\
    &S:=S^w_\delta\cup S^p_\delta.%\\
    %&u:=u^w\cup u^p,\\
    %&c:=c^w\cup c^p.
\end{align*} 
Enumerate the set $\{\alpha\}\cup\{\rho_\delta(\alpha'):\alpha'\in\dom(f^p_\delta)-\alpha\}\cup\{\omega_1\}$ in a strictly increasing order as
\[
    \alpha=\rho_0<\rho_1<\dots<\rho_n=\omega_1.
\]
Let $X_0:=\dom(f)\cap M^\delta_\alpha=\dom(f^w_\delta)$ and for $k\in[1,n)$, let
    \[
        X_k:=\{\alpha'\in\dom(f):f(\alpha')\in M^\delta_{\rho_{k+1}}-M^\delta_{\rho_k}\}.
    \]
Note that:
    \begin{itemize}
        \item $X_k\subseteq M^\delta_{\rho_k}$,
        \item $f[X_k]\subseteq M^\delta_{\rho_{k+1}}- M^\delta_{\rho_k}$,
        \item every element in $f[X_k]$ has label at most $\rho_k$,
        \item $\dom(f)=X_0\cup\dots\cup X_n$,
    \end{itemize}
As said above, the goal is to find a common extension $v\leq w\rest\delta,p\rest\delta$ such that $v$ decides edges in the set $\dom(f)=X_0\cup\dots\cup X_n$, and $v(0)$ decides edges in $u$ exactly in such a manner that makes $f:\dom(f^w_\delta)\cup\dom(f^p_\delta)\to\dot\GG$ edge-preserving. 
If such $v$ can be found, we can let $q:=v^\smallfrown(f,S^w_\delta\cup S^p_\delta)$, and this will be a condition and a common extension of $w$ and $p$. We find such a $v$ in $n$ many steps.

\vv 

We now define conditions $v_k\in\P_\delta\cap M^\delta_{\rho_k}$ by recursion on $k<n$. Let $v_0:=w\rest\delta$. Assume that $v_k$ was defined and satisfies:
\begin{enumerate}
	\item $v_k\in\P_\delta\cap M^\delta_{\rho_k}$,
	\item $v_k$ decides edges in the set $\dom(f^w_\delta)\cup X_0\cup\dots\cup X_{k-1}$, and for $\alpha_1,\alpha_2\in\dom(f^w_\delta)\cup X_0\cup\dots\cup X_{k-1}$,
    \begin{itemize}
    	\item $ v_k\Vdash E^{\dot H_\delta}(\alpha_1,\alpha_2)\quad\implies\quad E^{v_k}(\{f(\alpha_1),f(\alpha_2)\})=1$,
	\item $v_k\Vdash \neg E^{\dot H_\delta}(\alpha_1,\alpha_2)\quad\implies\quad E^{v_k}(\{f(\alpha_1),f(\alpha_2)\})=0$
\end{itemize}
\end{enumerate}

Up to extending $v_k$ in the model $M^\delta_{\rho_k}$, we may assume that it also decides edges in the set $\dom(f^w_\delta)\cup X_k$.

We proceed in two steps: we first extend $[p\rest\delta]^\delta_{\rho_{k+1}}$ to some $\tilde p$ by adding some edges to $E^p$, and then we obtain $v_{k+1}$ by amalgamating $v_k$ with $\tilde p$ in the model $M^\delta_{\rho_{k+1}}$.

To this end, let $\tilde E$ be the minimal extension of $E^p\cap M^\delta_{\rho_{k+1}}$  that satisfies that for every $\alpha_1\in\dom(f^w_\delta)$ and $\alpha_2\in X_k$:
\begin{itemize}
	\item if $v_k\Vdash ``E^{\dot H_\delta}(\alpha_1,\alpha_2)$, then $\tilde E(\{f^w_\delta(\alpha_1),f^p_\delta(\alpha_2)\})=1$,
	\item if $v_k\Vdash ``\neg E^{\dot H_\delta}(\alpha_1,\alpha_2)$, then $\tilde E(\{f^w_\delta(\alpha_1),f^p_\delta(\alpha_2)\})=0$,
\end{itemize}
Note that $f^p_\delta(\alpha_2)\in M^\delta_{\rho_{k+1}}-M^\delta_{\rho_k}$ and furthermore, if $\alpha_1\in\dom(f^w_\delta)-\dom(f^p_\delta)$, then the pair $\{f^w_\delta(\alpha_1),f^p_\delta(\alpha_2)\}$ is not yet in $\dom(E^p)$, by the assumption that $p$ is super-nice with respect to $M^\delta_\alpha$. Let $\tilde p$ be the condition obtained from $[p\rest\delta]^\delta_{\rho_{k+1}}$ by replacing $E^p\cap M^\delta_{\rho_{k+1}}$ by $\tilde E$.

\begin{claim}\label{claim:sp:1} $\tilde p$ is super-nice with respect to $M^\delta_{\rho_k}$ and $[\tilde p]^\delta_{\rho_k}=[p\rest\delta]^\delta_{\rho_k}$.
\end{claim}
\begin{proof}[Proof of Claim \ref{claim:sp:1}]
	It is clear that $[\tilde p]^\delta_{\rho_k}=[p\rest\delta]^\delta_{\rho_k}$, because we only added edges where the other end lies in $M^\delta_{\rho_{k+1}}-M^\delta_{\rho_k}$. We show that $\tilde p$ is super-nice with respect to $M^\delta_{\rho_k}$. Suppose not. By assumption $[p\rest\delta]^\delta_{\rho_{k+1}}$ is super-nice with respect to $M^\delta_{\rho_k}$, so by minimality of $\tilde E$, there is a path from $M^\delta_{\rho_k}$ to some $M^\gamma_\beta$ and there are $\xi\in \gamma\cap M^\gamma_\beta$ and $\alpha'\in\dom(f^p_\gamma)-\beta$ such that $f^p_\xi(\alpha')\in f^p_\delta[X_k]$. But every node in the set $f^p_\delta[X_k]$ must have label at most $\rho_k$, and only nodes in $\rho_k$ can be mapped to nodes with label $\rho_k$. But $\alpha'\not<\rho_k$, because $\alpha'\in\dom(f^p_\delta)-\beta$ and $\beta\geq\rho_k$. Thus $\alpha'$ cannot be mapped to a node in $f^p_\delta[X_k]$. Hence $\tilde p$ must be super-nice with respect to $M^\delta_{\rho_k}$.
	
\end{proof}

Now  $[\tilde p]^\delta_{\rho_k}$ is a residue of $\tilde p$ into $M^\delta_{\rho_k}$ and $v_k$ extends it. Thus there is a condition 
\[
v_{k+1}\leq v_k,\tilde p
\]
in $\P_\delta$. By elementarity, we may assume $v_{k+1}\in\P_\delta\cap M^\delta_{\rho_{k+1}}$. This ends Step $k+1$.

\vv 

Finally, define
\[
q:={v_n}^\smallfrown(f,S).
\]
Then $q$ is a condition in $\P_{\delta+1}$ and a common extension of $w$ and $p$. This ends the proof of Proposition \ref{prop:sp}.

\end{proof}

\begin{cor}
    For every $\delta<\omega_2$, the poset $\P_\delta$ is strongly proper and has $ccc$. The poset $\P_{\omega_2}$ has ccc.
\end{cor}
\begin{proof}
    Let $\lambda>\theta$ be a regular cardinal and let $M\elem H(\lambda)$ be countable model such that $\P_\delta\in M$ and $\alpha:=\omega_1\cap M\in\EE_\delta$. It suffices to show that $1_{\P_\delta}$ is strongly $(\P_\delta, M)$-generic. We have
    \[
     M\cap\P_\delta= M^\delta_\alpha\cap\P_\delta.
    \]
    But since $\P_\delta$ is strongly proper with respect to the model $ M^\delta_\alpha$, this implies the conclusion.
\end{proof}

\begin{cor}
    $2^\omega=\aleph_2$ in $V^{\P_{\omega_2}}$.
\end{cor}
\begin{proof}
    Each $\P_\delta$, $\delta<\omega_2$, is strongly proper and thus adds reals. The direction $2^\omega\leq\aleph_2$ follows from the fact that $\P_{\omega_2}$ has size $\aleph_2$.
\end{proof}

\begin{cor}
    In $V^{\P_{\omega_2}}$:
    \begin{enumerate}
        \item $\dot\GG$ is a universal countably chromatic graph of size $\aleph_1$,
        \item $2^\omega=\aleph_2$.
    \end{enumerate}
\end{cor}

\noindent We have proved:

\begin{thm}\label{thm}
    It is consistent that the class of countaby chromatic graphs of size $\aleph_1$ has a universal object and $2^{\omega}=\aleph_2$.
\end{thm}

\bibliographystyle{plain}
\bibliography{bib}

\end{document}